\documentclass[preprint,superscriptaddress,aps,pra,superscriptaddress]{revtex4-2}
\usepackage{natbib}
\usepackage{amsmath}
\usepackage{amssymb}
\usepackage{amsthm}
\usepackage{nicefrac}
\usepackage{graphicx}
\usepackage{bm}
\usepackage{url}
\DeclareMathOperator{\sdot}{\scriptscriptstyle\odot}
\newcommand{\bmm}[1]{\bm{\mathrm{#1}}}
\begin{document}
\title{On Sion's Minimax Theorem, Compact QCQPs, and Wave Scattering Optimization}
\author{Sean Molesky}
\email{sean.molesky@polymtl.ca}
\affiliation{Department of Engineering Physics, Polytechnique 
Montr\'{e}al, Montr\'{e}al, Qu\'{e}bec H3T 1J4, Canada}
\author{Pengning Chao}
\author{Alejandro W. Rodriguez}
\affiliation{Department of Electrical and Computer Engineering, Princeton 
University, Princeton, New Jersey 08544, USA}
\begin{abstract}
  \noindent
  In these notes, we examine certain implications of Sion's minimax theorem for compact quadratically constrained quadratic programs (QCQPs), particularly QCQPs arising in the context of optimizing wave scattering, in relation to Lagrangian duality. 
  The discussion puts forward an alternative ``dual'' understanding of optimization for wave phenomena that anticipates the realization of algorithmic (inverse design) methods attaining a guaranteed degree of global optimality for common figures of merit appearing in applied photonics, acoustics, and quantum mechanics.
\end{abstract}
\maketitle 
Building out of underlying relations to power transfer, information transfer, and the principle of stationary action~\cite{landau1976mechanics}, many of the common figures of merit considered in applied wave physics are described by quadratic functions of the (induced) response fields, Eq.~\eqref{quadFunc} Refs.~\cite{lu2012objective,miller2013complicated,roes2012acoustic,simon2015quantum,ganahl2017continuous}. 
When viewed through the lens of optimization theory~\cite{angeris2021convex,chao2021physical}, noting that the relations of scattering physics may also be simply phrased as quadratic constraint functions when the potential is linear, it is thus possible to state many of the optimization problems that appear in the design of wave devices as quadratically constrained quadratic programs (QCQPs)---c.f. Refs.~\cite{angeris2019computational,gustafsson2020upper,molesky2020t} and the topical reviews presented in Refs.~\cite{angeris2021heuristic,chao2021physical,liska2021fundamental}. 

Because the complexity class of QCQPs is $\mathsf{NP}$-hard~\cite{aaronson2003p,park2017general}, and the range of QCQPs that can be reformulated as the optimization of some abstract ``linear scattering potential'' is quite broad~\cite{angeris2021heuristic,chao2021physical}, this connection between QCQPs and wave device design would seem to indicate that, for certain problem classes, determination of a globally optimal solution (e.g. a best geometry for an electromagnetic device) becomes computationally intractable for even a modest number of degrees of freedom---although new approaches may yield substantial performance improvements~\cite{jensen2011topology,molesky2018inverse,christiansen2021inverse,digani2022framework,zhao2022efficient}, the prospect of universally determining globally optimal solutions is dubious. 
However, in opposition to the inferences one may sensibly draw from this likelihood, recent investigations have found that standard duality and semi-definite heuristics reliably give very good, often exact, approximations to QCQPs arising in the context of determining limits on achievable electromagnetic scattering~\cite{miller2019waves,kuang2020maximal,schab2020trade,molesky2020hierarchical,kuang2020computational,jelinek2020sub,capek2021fundamental}.
\begin{figure}[t!]
 \vspace{-12pt}
 \centering
 \includegraphics[width=1.0\columnwidth]{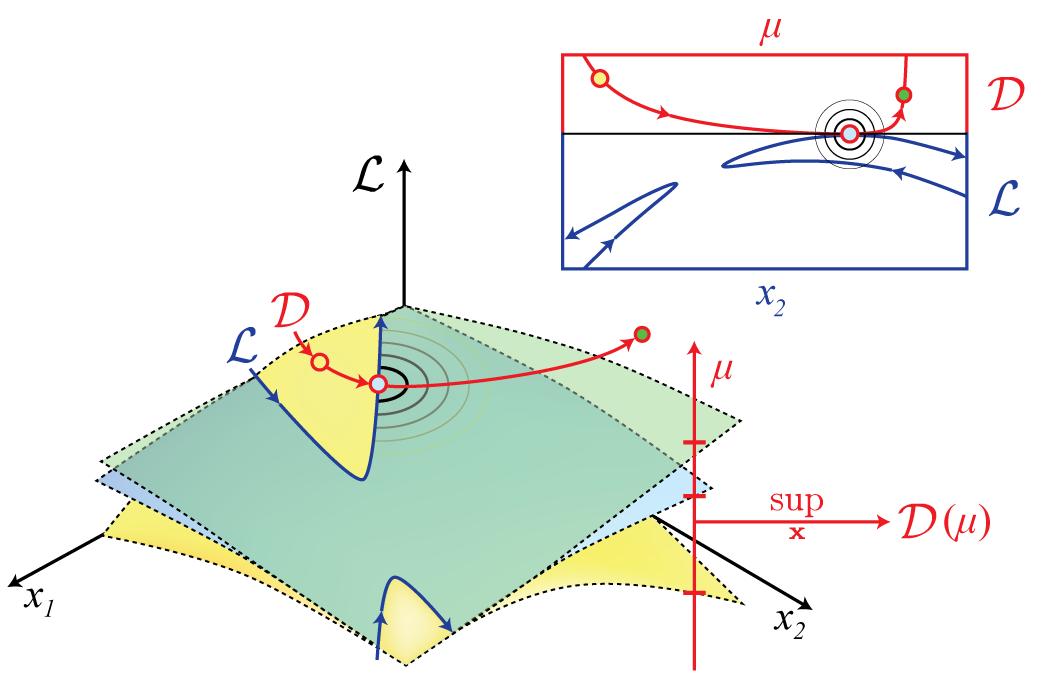}
 \caption{\bmm{Pictorial Lagrangian duality.} 
  The left part of the figure depicts three sections ($\mu$ multiplier values) of a schematic optimization Lagrangian $\mathcal{L}$ for a two dimensional $\left(x_{1},x_{2}\right)$ QCQP subject to a single feasible equality constraint, the dark blue path intersecting the three sections. 
  The value of the (convex) unconstrained dual $\mathcal{D}$ within each section is determined by the supremum of $\mathcal{L}$ over the associated manifold (represented by the colored balls) is necessarily at least as large as the maximum of $\mathcal{L}$ over the constraint set (represented by the ripple circles). 
  In the illustration, as indicated by the upper inset, the infimum of $\mathcal{D}$ (the supremum of the light blue section) is exactly equal to the maximum of $\mathcal{L}$ along the constraint path.
  For a general QCQP no such agreement between the maximum of $\mathcal{L}$ within the set of points satisfying all constraints and the infimum of $\mathcal{D}$ need occur, with the two values generally differing by a ``duality gap''. 
  An optimization program exhibits ``strong duality'' when the two values agree.
  }
  \label{dualFig}
% \vspace{-24pt}
\end{figure}

Here, we submit an intuitive explanation as to why these heuristics should usually work well for QCQPs with at least one compact constraint, particularly QCQPs originating from wave equations, and establish several relevant connections for Lagrangian dual programs in relation to such ``compact'' QCQPs. 
Making a minor generalization of the frameworks of Refs.~\cite{molesky2020hierarchical,kuang2020computational}, we first demonstrate that the optimization Lagrangian of any compact QCQP, covering a range of practical use cases in photonics, acoustics, and quantum mechanics~\cite{molesky2020t,zhang2021conservationlawbased,angeris2021convex,chao2021physical}, is consistent with the assumptions of Sion's minimax theorem---the order of maximization and minimization in the dual is freely interchangeable. 
Via this equivalence, we then find a simple condition related to the position of the dual solution relative to a boundary, the boundary set by the conservation of resistive power for wave scattering optimizations~\cite{chao2021physical}, under which strongly duality holds, Fig.~\ref{dualFig}. 
Expanding off this complementary understanding of Lagrangian dual programs, we then describe how strong duality could be retroactively achieved in any compact, feasible, QCQP by modifying the linear part of the objective function. 
That is, if the objective is allowed to vary, then a globally optimal solution in fact can be found. 
Beyond theoretical interest, the alternative perspective of optimization design forwarded herein may be conceivably leveraged to create algorithmic methods achieving a certifiable degree of optimality using only gradient information~\cite{angeris2019computational,chao2021physical,taylor2021optimal}. 
\subsection*{Notation and Physical Correspondence}
Throughout the text $\Omega = \mathbb{C}^{d}$ with $d\in\mathbb{N}$.
Bold lowercase letters represent either vectors (columns of complex numbers) in $\Omega$, when undecorated, or complex valued linear functionals (conjugated rows of complex numbers) on $\Omega$, when proceeded by a $\dagger$ superscript.  
Bold capital letters are used for linear maps (matrices) from $\Omega$ to itself, which are ``static'' when affixed with a Latin letter subscript and ``variable'' when affixed with a Greek letter subscript. 
$\bmm{I}$ is used to denote the identity operator on $\Omega$, and $\bmm{Z}$ is the constant map taking $\Omega$ to the zero vector. 
Standard matrix operations are implied by juxtaposition. 
Unless otherwise stated, all other symbols follow usual conventions: $\bmm{A}\succ \epsilon$ (resp. $\bmm{A}\succeq \epsilon$) means that $\bmm{A} - \epsilon \bmm{I}$ is Hermitian and positive definite (resp. semi-definite); $\ni$ reads as ``such that'', $\forall$ as ``for all'', $\backslash$ as ``removing'', $|$ as ``restricted to'', $\Rightarrow$ as ``implies that'', $\land$ as ``and'', and $\in$ as ``is an element of''; $\Re\left(z\right)$ and $\Im\left(z\right)$ are the real and imaginary parts, respectively, of the complex number $z$. 
In stating that a function $f\left(\bmm{t}\right)$ is quadratic, $f:\Omega\rightarrow\mathbb{R}$, we mean that it may be written as 
\begin{equation}
  f\left(\bmm{t}\right) = 2~\Re\left(\bmm{t}^{\dagger}\bmm{s}\right) - \bmm{t}^{\dagger}\bmm{A}\bmm{t} + v
  \label{quadFunc}
\end{equation}
with $\bmm{s}$, $\bmm{A}$ and $v\in\mathbb{R}$ known, and $\bmm{A}$ Hermitian. 
If $\bmm{A}\succ 0$ (resp. $\bmm{A}\succeq 0$), then we will say that $f_{-}\left(\bmm{t}\right)$ is a positive definite (resp. semi-definite) quadratic function. 
Relatedly, if $\bmm{A}\succ \epsilon$ (resp. $\bmm{A}\succeq\epsilon$), we will say that $f\left(\bmm{t}\right)$ is $\epsilon$-positive definite (resp. semi-definite). 
$\mathsf{h}$ and $\mathsf{s}$ superscripts are used for the Hermitian and skew-Hermitian parts of a linear map; i.e., $\bmm{A}^{\mathsf{h}} = \left(\bmm{A} + \bmm{A}^{\dagger}\right)/2$ and $\bmm{A}^{\mathsf{s}} = \left(\bmm{A} - \bmm{A}^{\dagger}\right)/2$, with $\!^{\dagger}$ the conjugate transpose of the antecedent operator. 
$\circledast$ superscripts denote optimal solutions based on context, and a $\partial$ superscript denotes the boundary of a set. 
$\lVert -\rVert_{_\mathtt{E}}$ is used for the (usual) Euclidean norm and $\lVert -\rVert_{_\mathtt{O}}$ for the operator norm. 

In regards to interpretation, the following correspondences can be made between the symbols appearing below and the standard quantities of electrodynamics. 
Up to unit conversion factors, $\bmm{w}$ can be equated with a total electric field; $\bmm{s}$ with an incident electric field; $\bmm{M}$ with the d'Alembert operator; $\bmm{V}_{\bm{\rho}}$ with a variable electrical susceptibility profile, which would typically define a material  geometry; \bmm{T} and $\bmm{G}$ are the T-operator and Green's function~\cite{tsang2004scattering}; and $\bmm{t}_{\rho}$, again up to unit conversion factors, is the total electric polarization current density.  

\subsection*{Assumptions}
The subsequent assumptions, developed from Refs.~\cite{angeris2019computational,molesky2020hierarchical,gustafsson2020upper,kuang2020maximal,trivedi2020bounds,molesky2020hierarchical,kuang2020computational,schab2020trade,jelinek2020sub} and summarized in Refs.~\cite{angeris2021heuristic,chao2021physical}, are made in most statements appearing hereafter. 
\\ \\
\emph{Linear response}---In connection to physical device design, the optimization problem is described by some quadratic objective function $f_{\circ}\left(\bmm{w}\right)$ of a total response field $\bmm{w}$, with $v_{\circ} = 0$, subject to a system of linear constraints $\bmm{M}\left(\bmm{w} + \bmm{s}\right) = \bmm{V}_{\rho}\left(\bmm{w} + \bmm{s}\right)$, where $\bmm{V}_{\rho}$ is a variable (designable) scattering potential, and $\bmm{s}$ is a known initial field:
\begin{align}
  &\max_{\bmm{V}_{\rho}}~f_{\circ}\left(\bmm{w}\right)\ni \nonumber \\
  &\bmm{M}\left(\bmm{w} + \bmm{s}\right) = \bmm{V}_{\rho}\left(\bmm{w} + \bmm{s}\right).
  \label{basicLinOpt}
\end{align}
\\ 
\emph{Comments}---The assumptions of linear response, potential separability, and passivity, are only important in so far as justifying the pervading assumption of feasibility, and the general applicability of Eq.~\eqref{genCon2} to scattering optimization problems. 
The results stated under \emph{Dual Properties for Compact QCQPs}, \emph{Criteria for Strong Duality in Compact QCQPs}, and \emph{``Dual'' Properties under Minimax Exchange} can be adapted to any feasible QCQP wherein at least one of the imposed constraints is positive definite without any additional technical criteria. 
\\ \\
\emph{Potential separability}---There is a collection $\left\{\bmm{D}_{j}\right\}_{j\in\mathtt{J}}$ of ``design'' matrices respecting the following properties. 
For all $j\in \mathtt{J}$, $\bmm{D}_{j}$ is restricted to a subspace $\mathtt{d}_{j}\subseteq \Omega$ and the set $\left\{\mathtt{d}_{j}\right\}_{j\in\mathtt{J}}$ partitions $\Omega$: if $j\neq k$, $\mathtt{d}_{k}\cap \mathtt{d}_{j} = \emptyset$ and $\bigcup_{j\in J}~\mathtt{d}_{j} = \Omega$. 
To every scattering potential $\bmm{V}_{\rho}$ there is a ``parameterization'' $\bm{\rho} \in \left\{0,1\right\}^{\mathtt{J}}$ such that 
\begin{equation}
  \bmm{V}_{\rho}^{-1} = 
  \sum_{j\in \mathtt{J}} 
  \rho_{j} \bmm{D}_{j}, 
  \label{svdDesign}
\end{equation}
with $\bmm{V}_{\rho}^{-1}$ denoting the pseudo-inverse of $\bmm{V}_{\rho}$.
\\ \\
\emph{Comments}---Potential separability is a property of any scattering theory where the potential may only take on a single nonzero value and does not mix distinct spatial points (although much of what is done below can be extended to broader classes of potentials). 
If the potential is ``local'' in the sense that the polarization generated at any spatial point requires only knowledge of the total vector field at that point, and whether the potential is absent or present at the same point, then the scattering potential is separable as defined above. 
In the discussion presented below, some subcollection of the $\bm{\rho}$ may also be considered as fixed, i.e. $\rho_{-} = 0$ or $\rho_{-} = 1$. 
The assumption that $\bigcup_{j\in J}~\mathtt{d}_{j} = \Omega$ may then be relaxed in order to give $\Omega$ a simpler description, or to treat problems where part of the potential is given and unalterable. 
\\ \\
\emph{Passivity}---There is a positive real number $\epsilon$ such that $i\bmm{D}^{\mathsf{s}}_{j}|_{\mathtt{d}_{j}}\succ \epsilon$ for every $j\in\mathtt{J}$.
\\ \\
\emph{Comments}---$-i\bmm{R}^{\mathsf{s}}$, for any physical response $\bmm{R}$ (a susceptibility, a Green's function, etc.), describes the power lost in generating the response field to ``external'' mechanisms~\cite{landau2013statistical}, e.g. in electrodynamics the absorption of radiation or power escaping beyond a simulation boundary. 
Hence, whenever $\left\{\bmm{D}_{j}\right\}_{j\in\mathtt{J}}$ is used to represent the admittance associated with a collection of causal and passive scattering potentials, as would happen in the conception of a passive wave scattering device, the passivity criterion should hold. 
$i\bmm{D}^{\mathsf{s}}_{j}|_{\mathtt{d}_{j}}\succ \epsilon$ asserts that in creating the response field a certain amount of power, in fixed proportion to the square of the response field amplitude, must be irrecoverably lost to the ambient environment~\cite{molesky2020t,molesky2020hierarchical}. 
For general QCQPs, the assumption of passivity is replaced with the assumption that one of the imposed quadratic constraints is positive definite. 
\\ \\
\emph{Feasibility}---Feasible points exist; i.e., there are fields (points) that satisfy all imposed constraints. 
\\ \\
\emph{Comments}---As it is not possible to construct constraints from a scattering theory that exclude the response fields of physically realizable potentials, the assumption of feasibility holds in any physical optimization problem. 
Trivially, if $\bmm{V}_{\rho} = \bmm{0}$ (the no scattering potential) is admissible, then the zero vector is feasible. 
The assumption of feasibility is a far stronger restriction on general compact QCQPs. 
\subsection*{Quadratic Constraints for Wave Scattering}
Leveling the standard perspective of beginning with empty space and then adding material~\cite{felsen1994radiation,tsang2004scattering,newton2013scattering}, set $\bmm{V}_{\bm{b}}$ to be the potential of some specific design, and take $\tilde{-}$ quantities to denote scattering relative to $\bmm{V}_{\bm{b}}$. 
Letting alike undecorated symbols mark the special case of the empty (``all-off'') potential, take
\begin{align}
  \bmm{t} = \bmm{T}\bmm{s}
  ~\land~
  \tilde{\bmm{t}} = \tilde{\bmm{T}}\tilde{\bmm{s}},
\end{align}
with $\bmm{s}$ a given incident field, $\tilde{\bmm{s}} = \left(\bmm{I} -\bmm{G}\bmm{V}_{\bm{b}}\right)^{-1}\bmm{s} = \bmm{X}_{\bm{b}}\bmm{s}$, and $\bmm{T}$ and $\tilde{\bmm{T}}$ defined by the $\bmm{T}$-operator relations
\begin{align}
  \bmm{I}_{\bm{\rho}} = \bmm{I}_{\bm{\rho}}\left(\bmm{V}_{\bm{\rho}}^{-1} - \bmm{G}\right)\bmm{T}
  ~\land~
  \bmm{I}_{\bm{b\rho}} = \bmm{I}_{\bm{b\rho}}\left(\tilde{\bmm{V}}^{-1}_{\bm{\rho}} - \tilde{\bmm{G}}\right)\tilde{\bmm{T}},
  \label{tBasic}
\end{align}
where $\bmm{I}_{\bm{\rho}}$ is the characteristic operator of the design $\bm{\rho}$, and $\bmm{I}_{\bm{b\rho}}$ is similarly defined for the design $\bm{\rho}$ relative to the background $\bm{b}$. 
Because both $\bmm{V}_{\bm{\rho}}$ and $\tilde{\bmm{V}}_{\bm{\rho}}$ respect---are block diagonal in---the decomposition of $\Omega$ given by $\left\{\mathtt{d}_{j}\right\}_{j\in\mathtt{J}}$, the specification of the design potential in both $\bmm{I}_{\bm{\rho}}$ and $\bmm{V}_{\bm{\rho}}^{-1}$ (resp. $\bmm{I}_{\bm{b\rho}}$ and $\tilde{\bmm{V}}^{-1}_{\bm{\rho}}$) is redundant, and the forms given in Eq.~\eqref{tBasic}, acting from the left with $\bmm{T}$ and $\tilde{\bmm{T}}$ respectively, are equivalent to 
\begin{align}
  \bmm{T}^{\dagger}\bmm{P} &= \bmm{T}^{\dagger}\bmm{P}\left(\bmm{V}^{-1} - \bmm{G}\right)\bmm{T} = \bmm{T}^{\dagger}\bmm{P}\bmm{U}\bmm{T}~\land
  \nonumber \\
  \tilde{\bmm{T}}^{\dagger}\bmm{P} &= \tilde{\bmm{T}}^{\dagger}\bmm{P}\left(\tilde{\bmm{V}}^{-1} - \tilde{\bmm{G}}\right)\tilde{\bmm{T}} = \tilde{\bmm{T}}^{\dagger}\bmm{P}\tilde{\bmm{U}}\tilde{\bmm{T}}, 
  \label{tGen}
\end{align}
where $\bmm{V}^{-1}$ is the inverse of the scattering potential for the ``all on'' design, $\tilde{\bmm{V}} = \bmm{V} - 2\bmm{V}_{\bm{b}}$, $\bmm{P}$ is any map that respects the decomposition of $\Omega$ given by $\left\{\mathtt{d}_{j}\right\}_{j\in\mathtt{J}}$---any map that does not mix elements from distinctly indexed $\mathtt{d}_{j}$ subspaces, and $\bmm{U}$ and $\tilde{\bmm{U}}$ are defined through Eq.~\eqref{tGen}.
As the total field $\bmm{w}$ linked with any particular design can not change depending on how the background is described, the definitions of $\bmm{T}$, $\bmm{G}$, $\tilde{\bmm{T}}$, and $\tilde{\bmm{G}}$ enforce that 
\begin{equation}
  \bmm{s} + \bmm{G}\bmm{t} = 
  \tilde{\bmm{s}} + \tilde{\bmm{G}}\tilde{\bmm{t}}.
  \label{totalFieldEquivalence}
\end{equation}
Correspondingly, noting that
\begin{equation}
  \tilde{\bmm{G}}^{-1} = \bmm{G}^{-1} - \bmm{V}_{\bm{b}},
  \label{gInvEq}
\end{equation}
the definition of $\bmm{W}_{\bm{b}} = \left(\bmm{I} - \bmm{V}_{\bm{b}}\bmm{G}\right)^{-1} =\bmm{G}^{-1}\bmm{X}_{\bm{b}}\bmm{G}$, 
\begin{equation}
  \tilde{\bmm{G}} = \bmm{G}\bmm{W}_{\bm{b}} \Rightarrow\bmm{G}^{-1}\tilde{\bmm{G}}=\bmm{W}_{\bm{b}},
  \label{wImplications}
\end{equation}
implies that Eq.~\eqref{totalFieldEquivalence} may be rewritten as 
\begin{align}
  \tilde{\bmm{t}} &= \tilde{\bmm{G}}^{-1}\left(\bmm{s}-\tilde{\bmm{s}}\right) + 
  \tilde{\bmm{G}}^{-1}\bmm{G}\bmm{t}
  \nonumber \\
  &= \tilde{\bmm{G}}^{-1}\bmm{s}-\tilde{\bmm{G}}^{-1}\bmm{X}_{\bm{b}}\bmm{s} + 
  \bmm{W}_{\bm{b}}^{-1}\bmm{t}
  \nonumber \\
  &= \tilde{\bmm{G}}^{-1}\bmm{s}-\tilde{\bmm{G}}^{-1}\bmm{G}\bmm{W}_{\bm{b}}\bmm{G}^{-1}\bmm{s} + 
  \bmm{W}_{\bm{b}}^{-1}\bmm{t}
  \nonumber \\
  &= \tilde{\bmm{G}}^{-1}\bmm{s}-\bmm{G}^{-1}\bmm{s} + 
  \bmm{W}_{\bm{b}}^{-1}\bmm{t} = -\bmm{V}_{b}\bmm{s} + \bmm{W}_{\bm{b}}^{-1}\bmm{t}.
  \label{tTransform}
\end{align}
Encapsulating the second form of Eq.~\eqref{tGen} with $\tilde{\bmm{s}}^{\dagger}\dots\tilde{\bmm{s}}$, use of Eq.~\eqref{totalFieldEquivalence} then gives
\begin{equation}
  \tilde{\bmm{t}}^{\dagger}\bmm{P}\left(\bmm{s} + \bmm{G}\bmm{t}\right) 
  = \tilde{\bmm{t}}^{\dagger}\bmm{P}\tilde{\bmm{V}}^{-1}\tilde{\bmm{t}}. 
  \label{genCon1}
\end{equation}
Through Eq.~\eqref{tTransform}, the $\bmm{t}$ independent part of Eq.~\eqref{genCon1} is 
\begin{equation}
  -\bmm{s}^{\dagger}\bmm{V}_{\bm{b}}^{\dagger}\bmm{P}\bmm{s} =
  \bmm{s}^{\dagger}\bmm{V}_{\bm{b}}^{\dagger}\bmm{P}\tilde{\bmm{V}}^{-1}
  \bmm{V}_{\bm{b}}\bmm{s},
  \label{constPart}
\end{equation}
which immediately cancels as either $\bmm{V}_{\bm{b}}|\mathtt{d}_{j} = \bmm{t}|\mathtt{d}_{j}$ or $-\tilde{\bmm{V}}|\mathtt{d}_{j} = \bmm{V}_{\bm{b}}|\mathtt{d}_{j}$. 
Analogously, the quadratic term of Eq.~\eqref{genCon1} after application of Eq.~\eqref{tTransform} is  
\begin{align}
  &\bmm{t}^{\dagger}\bmm{W}_{\bm{b}}^{-1\dagger}\bmm{P}
  \left(\tilde{\bmm{V}}^{-1}-\tilde{\bmm{V}}^{-1}\bmm{V}_{\bm{b}}\bmm{G} - \bmm{G}
  \right)\bmm{t} = 
  \nonumber \\
  &\bmm{t}^{\dagger}\bmm{W}_{\bm{b}}^{-1\dagger}\bmm{P}\tilde{\bmm{V}}^{-1}\bmm{W}_{\bm{c}}^{-1}\bmm{t}
\end{align}
where $\bmm{W}_{\bm{c}}^{-1} = \bmm{I} - \bmm{V}_{\bm{c}}\bmm{G}$, and $\bmm{V}_{\bm{c}}$ is ``complementary'' to $\bmm{V}_{\bm{b}}$: $\bmm{V}_{\bm{c}}|\mathtt{d}_{j} = \bmm{t}|\mathtt{d}_{j}$ whenever $\bmm{V}_{\bm{b}}|\mathtt{d}_{j} = \bmm{D}|\mathtt{d}_{j}$ and $\bmm{V}_{\bm{c}}|\mathtt{d}_{j} = \bmm{D}|\mathtt{d}_{j}$ whenever $\bmm{V}_{\bm{b}}|\mathtt{d}_{j} = \bmm{t}|\mathtt{d}_{j}$.
The terms of Eq.~\eqref{genCon1} in which only one of $\bmm{t}$ and $\bmm{t}^{\dagger}$ is present are
\begin{align}
  &\bmm{t}^{\dagger}\bmm{W}_{\bm{b}}^{-1\dagger}\bmm{P}\left(\bmm{I} +\tilde{\bmm{V}}^{-1}\bmm{V}_{\bm{b}}\right)\bmm{s} =
  \bmm{t}^{\dagger}\bmm{W}_{\bm{b}}^{-1\dagger}\bmm{P}\tilde{\bmm{V}}^{-1}\bmm{V}_{\bm{c}}\bmm{s}~\land
  \nonumber \\
  &\bmm{s}^{\dagger}\bmm{V}_{\bm{b}}^{\dagger}\bmm{P}\left(\tilde{\bmm{V}}^{-1}\bmm{W}_{\bm{b}}^{-1}-\bmm{G}\right)\bmm{t}
  =
  \bmm{s}^{\dagger}\bmm{V}_{\bm{b}}^{\dagger}\bmm{P}\tilde{\bmm{V}}^{-1}\bmm{W}_{\bm{c}}^{-1}\bmm{t}. 
\end{align}
Therefore, symmetrizing Eq.~\eqref{genCon1} and redefining $\bmm{P}$ as $\bmm{P}\tilde{\bmm{V}}^{-1}$ yields the constraint schema
\begin{align}
  &\Re\left[\bmm{t}^{\dagger}\left(\bmm{W}_{\bm{b}}^{-1\dagger}\bmm{P}\bmm{V}_{\bm{c}} + \bmm{W}_{\bm{c}}^{-1\dagger}\bmm{P}^{\dagger}\bmm{V}_{b}\right)\bmm{s}\right] 
  =
  \nonumber \\
  &\bmm{t}^{\dagger}\left(\bmm{W}_{\bm{b}}^{-1\dagger}\bmm{P}\bmm{W}_{\bm{c}}^{-1}\right)^{\mathsf{h}}\bmm{t}. 
  \label{genCon2}
\end{align}
In the case that the background is selected to be either empty or full Eq.~\eqref{genCon2} reverts trivially to the 
\begin{equation}
  \Re\left(\bmm{t}^{\dagger}\bmm{P}\bmm{s}\right) = \bmm{t}^{\dagger}\left(\bmm{P}\bmm{U}\right)^{\mathsf{h}}\bmm{t} 
  \label{prevGenConst}
\end{equation}
form studied in prior articles, Ref.~\cite{molesky2020hierarchical,kuang2020computational,jelinek2020sub,molesky2021comm}. 
For all other choices of background, however, the fact that $\bmm{W}_{\bm{b}}$ and $\bmm{W}_{\bm{c}}$ do not commute with the characteristic design maps $\bmm{I}_{\mathtt{d}_{j}}$ means that no simple transformation via $\bmm{P}$ can be used to directly equate Eq.~\eqref{genCon2} with Eq.~\eqref{prevGenConst}. 
Additionally, because imposing a constraint of the form given by Eq.~\eqref{genCon2} in no way alters the objective, no simple redefinition of $\bmm{t}$ can transform any QCQP subject to constraints of the form of Eq.~\eqref{genCon2} into a QCQP where all constraints are of the form of Eq.~\eqref{prevGenConst}. 

This is not to say that Eq.~\eqref{genCon2} contains additional physics compared to Eq.~\eqref{prevGenConst}. 
If every constraint admissible under Eq.~\eqref{prevGenConst} is used, any additional constraint will begin to definitively specify some designable aspect of the scattering potential~\cite{molesky2021comm}. 
Eq.~\eqref{genCon2} simply repackages Eq.~\eqref{prevGenConst}, as motivated by distinct (but ultimately equivalent) consistency requirements. 
Nevertheless, in situations where a variable but finite subset of constraints is imposed Eq.~\eqref{genCon2} supersedes Eq.~\eqref{prevGenConst}.
The added freedom offered by Eq.~\eqref{genCon2} may therefore have meaningful algorithmic consequences. 
\subsection*{Compact QCQPs for Wave Scattering}
Collating the above, making use of the scattering relation that up to unit conversion prefactors $\bmm{w} = \bmm{G}\bmm{t}$, take $f_{\mathsf{o}}\left(\bmm{t}\right)$ to be a quadratic objective function corresponding to some desirable physical property. 
Let $\mathtt{K} = \bigcup_{\bmm{P},\bmm{V}_{\bm{b}}}\left\{f_{p}\left(\bmm{t}\right)\right\}$ be a collection of quadratic equality constraints defined through Eq.~\eqref{genCon2} by distinct choices of $\bmm{P}$ and $\bmm{V}_{\bm{b}}$, and set $f_{\sdot}\left(\bmm{t}\right) = 2\Re\left(\bmm{t}^{\dagger}\bmm{s}_{\sdot}\right) - \bmm{t}^{\dagger}\bmm{A}_{\sdot}\bmm{t}$, with $\bmm{s}_{\sdot} = i\bmm{s}/2$ and $\bmm{A}_{\sdot} = \left(i\bmm{U}\right)^{\mathsf{h}}$. 
By the assumption of passivity $f_{\sdot}\left(\bmm{t}\right)$ is a positive definite quadratic function, so that the set where $f_{\sdot}\left(\bmm{t}\right)\geq 0$ is compact~\cite{bredon2013topology}. 
Taking $\tilde{\mathtt{K}} = \mathtt{K}\bigcup\left\{f_{\sdot}\left(\bmm{t}\right) = 0\right\}$, define $\mathcal{Q}\left(\Omega,f_{\circ},\tilde{\mathtt{K}}\right)$ as
\begin{align}
  &\max_{\bmm{t}\in\Omega}~f_{\mathsf{o}}\left(\bmm{t}\right)\ni~
  \nonumber \\
  &\left(\forall f_{p}\in\mathtt{K}\right)~f_{p}\left(\bmm{t}\right) = 0,~\land
  ~f_{\sdot}\left(\bmm{t}\right) = 0,
  \label{waveQCQP}
\end{align}
the QCQP associated with optimizing $f_{\circ}\left(\bmm{w}\right)$ subject to the constraint set $\tilde{\mathtt{K}}$.
Adopting $\bmm{s}_{\bm{\phi}} = \bmm{s}_{\mathsf{o}} + \phi_{\sdot}\bmm{s}_{\sdot} + \sum_{\mathtt{K}}\phi_{k} \bmm{s}_{k}$, $\bmm{A}_{\bm{\phi}} = \bmm{A}_{\mathsf{o}} + \phi_{\sdot}\bmm{A}_{\sdot} + \sum_{\mathtt{K}}\phi_{k} \bmm{A}_{k}$, and $v_{\bm{\phi}} = \phi_{\sdot}v_{\sdot} + \sum_{\mathtt{K}}\phi_{k}v_{k}$ 
\begin{align}
  \mathcal{L}\left(\bm{\phi}, \bmm{t}\right) =
  2~\Re\left(\bmm{t}^{\dagger}\bmm{s}_{\bm{\phi}}\right) - 
  \bmm{t}^{\dagger}\bmm{A}_{\bm{\phi}}\bmm{t} + v_{\bm{\phi}}
 \label{Lform}
\end{align}
is the Lagrangian of Eq.~\eqref{waveQCQP}. 
The convex Wolfe (or Lagrangian) dual program of Eq.~\eqref{waveQCQP}, Ref.~\cite{boyd2004convex}, is 
\begin{equation}
  \mathcal{W}\left(\Omega,f_{\circ},\tilde{\mathtt{K}}\right) 
  = \inf_{\bm{\phi}\in\Phi_{\epsilon}}\max_{\bmm{t}\in\mathtt{C}}\mathcal{L}\left(\bm{\phi},\bmm{t}\right) 
  = \inf_{\bm{\phi}\in\Phi_{\epsilon}}\mathcal{D}\left(\bm{\phi}\right), 
  \label{dualQCQP}
\end{equation} 
where $\mathcal{D}\left(\bm{\phi}\right) = \max_{\bmm{t}\in\mathtt{C}}\mathcal{L}\left(\bm{\phi},\bmm{t}\right)$, $\mathtt{C} = \left\{\bmm{t}\in\Omega~|~f_{\sdot}\left(\bmm{t}\right) \geq 0\right\}$, and $\Phi_{\epsilon} = \{\bm{\phi}\in\mathbb{R}^{\tilde{\mathtt{K}}}~|~\bmm{A}_{\bm{\phi}} \succeq \epsilon\}$ with $\epsilon \geq 0$. 

The Wolfe dual is conventionally defined as the infimum of $\sup_{\bmm{t}\in\Omega}\mathcal{L}$ over $\Phi_{0}$. 
However, the refinement to $\max_{\bmm{t}\in\mathtt{C}}$, and freedom offered by $\Phi_{\epsilon}$, is more convenient for our subsequent uses and makes no practical difference. 
$\bmm{A}_{\bm{\phi}}\succeq \epsilon$ is a more restrictive condition than $\bmm{A}_{\bm{\phi}}\succeq 0$ so $\inf_{\Phi_{\epsilon}}\sup_{\bmm{t}\in\Omega}\mathcal{L}$ remains a convex program bounding $\mathcal{Q}\left(\Omega,f_{\circ},\tilde{\mathtt{K}}\right)$. 
After making this first alteration, if $\bmm{t}^{\circledast}\not\in\mathtt{C}$ for some $\bm{\phi}$ then $f_{\sdot}\left(\bmm{t}^{\circledast}\right) < 0$, and by increasing the $\phi_{\sdot}$ component of $\bm{\phi}$ some sufficiently small amount the value of the dual may be reduced while remaining in $\Phi_{\epsilon}$.

In the subsequent four sections we will demonstrate that for any such compact QCQP (a) that a physically motivated boundary condition exist under which strong duality holds---that solving Eq.~\eqref{dualQCQP} solves Eq.~\eqref{waveQCQP} for many wave scattering optimization problems---and (b) that by retroactively modifying the linear $\bmm{s}_{\circ}$ part of the objective function $f_{\circ}\left(\bmm{t}\right)$ strong duality between Eq.~\eqref{waveQCQP} and Eq.~\eqref{dualQCQP} may be attained for any constraint set $\tilde{\mathtt{K}}$. 
\subsection*{Dual Properties for Compact QCQPs}
\noindent
\emph{Sion's minimax theorem~\cite{sion1958general,komiya1988elementary}.}
Let $\mathtt{C}$ be a compact convex subset of a topological vector space, $\Phi_{\epsilon}$ a convex subset of a topological vector space, and $\mathcal{L}$ a real valued function on $\Phi_{\epsilon}\times\mathtt{C}$. 
Suppose (a) that $\mathcal{L}\left(\cdot,\bmm{t}\right)$ is lower semicontinuous and quasi-convex on $\Phi_{\epsilon}$ for each $\bmm{t}\in\mathtt{C}$, and (b) that $\mathcal{L}\left(\bm{\phi},\cdot\right)$ is upper semicontinuous and quasi-concave on $\mathtt{C}$ for each $\bm{\phi}\in \Phi_{\epsilon}$. 
Then,
\begin{equation}
  \max_{\bmm{t}\in\mathtt{C}}\inf_{\bm{\phi}\in\Phi_{\epsilon}}\mathcal{L}\left(\bm{\phi},\bmm{t}\right) = 
  \inf_{\bm{\phi}\in\Phi_{\epsilon}}\max_{\bmm{t}\in\mathtt{C}}~\!\mathcal{L}\left(\bm{\phi},\bmm{t}\right).
  \label{sion} 
\end{equation}
\\ 
\emph{Lemma 1}---With $\mathtt{C}$ and $\Phi_{\epsilon}$ as stated under \emph{QCQPs for Wave Scattering}, Eq.~\eqref{sion} applies to compact QCQPs, i.e. QCQPs with one positive definite constraint: 
$$
  \mathcal{W}\left(\Omega,f_{\circ},\tilde{\mathtt{K}}\right) 
  = \inf_{\bm{\phi}\in\Phi_{\epsilon}}\mathcal{D}\left(\bm{\phi}\right) 
  = \max_{\bmm{t}\in\mathtt{C}}\inf_{\bm{\phi}\in\Phi_{\epsilon}}\mathcal{L}\left(\bm{\phi},\bmm{t}\right).
$$
\begin{proof}[Proof of lemma 1.]
  Since $\bmm{A}_{\sdot} = \left(i\bmm{U}\right)^{\mathsf{h}}\succ\epsilon$, $f_{\sdot}\left(\bmm{t}\right)$ is a strictly concave continuous function, Ref.~\cite{folland1999real}, that is bounded from above by $\bmm{s}^{\dagger}_{\sdot}\bmm{A}_{\sdot}^{-1}\bmm{s}_{\sdot}$.
  As such, $\mathtt{C}$ is compact and convex~\cite{bredon2013topology}.
  Suppose $\bmm{t},\bmm{r}\in\mathtt{C}$. 
  For any $s\in\left(0,1\right)$ 
  \begin{align}
    &sf_{\sdot}\left(\bmm{t}\right) + \left(1-s\right)f_{\sdot}\left(\bmm{r}\right) 
    =\nonumber \\ 
    &\Re\left[\left(s\bmm{t}+\left(1-s\right)\bmm{r}\right)^{\dagger}i\bmm{s}\right] - 
    s~\bmm{t}^{\dagger}\bmm{A}_{\sdot}\bmm{t} - \left(1-s\right)
    \bmm{r}^{\dagger}\bmm{A}_{\sdot}\bmm{r} =
    \nonumber \\
    &f_{\sdot}\left[s\bmm{t}+\left(1-s\right)\bmm{r}\right]
    -s\left(1-s\right)\left(\bmm{t}-\bmm{r}\right)
    \bmm{A}_{\sdot}\left(\bmm{t}-\bmm{r}\right) \Rightarrow
    \nonumber \\
    &f_{\sdot}\left[s\bmm{t}+\left(1-s\right)\bmm{r}\right]
    \geq sf_{\sdot}\left(\bmm{t}\right) + \left(1-s\right)f_{\sdot}\left(\bmm{r}\right)~\Rightarrow
    \nonumber \\
    &s\bmm{t}+\left(1-s\right)\bmm{r}\in\mathtt{C}.
    \nonumber
  \end{align}
  Similarly, $\Phi_{\epsilon}$ is convex since if $\bmm{A}_{\bm{\phi}}\succeq \epsilon$ and $\bmm{A}_{\bm{\psi}}\succeq \epsilon$ then $\left(\forall\bmm{t}\in\Omega~\land~\forall s\in\left(0,1\right)\right)$
  $$
    \bmm{t}^{\dagger}\left(s\left(\bmm{A}_{\bm{\phi}}-\epsilon\bmm{I}\right) +\left(1-s\right)\left(\bmm{A}_{\bm{\psi}}-\epsilon\bmm{I}\right)\right)\bmm{t} \geq 0
  $$ 
  $\Rightarrow s\bmm{A}_{\bm{\phi}} + \left(1-s\right)\bmm{A}_{\bm{\psi}} \in \Phi_{\epsilon}$.
  Because $\bmm{A}_{\bm{\phi}}\succeq \epsilon$, as shown above, $\mathcal{L}$ is concave and continuous function of $\bmm{t}$ for any fixed $\bmm{\phi}\in\Phi_{\epsilon}$. 
  For any fixed $\bmm{t}\in\mathtt{C}$, $\mathcal{L}$ is an affine function of $\bm{\phi}$, and is therefore convex and continuous. 
\end{proof}
Given that $\mathtt{C}$ is compact, the value of $\mathcal{D}\left(\bm{\phi}\right)$ is always finite, and there is always a point in $\mathtt{C}$ where the value of the supremum of the Lagrangian over $\bmm{t}\in\mathtt{C}$ (for any fixed $\bmm{\phi}$) is achieved, justifying the use of $\max$ in the above~\cite{bredon2013topology}. 
The properties of Eq.~\eqref{waveQCQP} also entail that $\mathcal{D}\left(\bm{\phi}\right)$ is a continuous, rather than semi-continuous, function on $\Phi_{\epsilon}$~\cite{rockafellar1970convex}, that the minimum of $\mathcal{D}\left(\bm{\phi}\right)$ on $\Phi_{\epsilon}$ (if it exists) is unique for all $\epsilon >0$~\cite{molesky2020hierarchical}, and that $\mathcal{D}\left(\bm{\phi}\right)$ is equivalently defined for any $\bm{\phi}\in\Phi_{\epsilon}$ as
\begin{equation}
 \mathcal{D}\left(\bm{\phi}\right) 
 = \bmm{t}^{\circledast\dagger}\bmm{A}_{\tilde{\bm{\phi}}}\bmm{t}^{\circledast},
 \label{dualEval}
\end{equation}
where $\bmm{A}_{\tilde{\bm{\phi}}}\succeq\epsilon$~\cite{beck2006strong}.
Here, $\bmm{t}^{\circledast}$ is defined by $\bmm{A}_{\tilde{\bm{\phi}}}\bmm{t}^{\circledast} = \textbf{y}_{\tilde{\bm{\phi}}}$, with $\bmm{t}^{\circledast}$ satisfying $f_{\sdot}\left(\bmm{t}^{\circledast}\right) = 0$ in the event that $\bmm{A}_{\tilde{\bm{\phi}}}$ is semi-definite.
In these definitions, $\tilde{\bm{\phi}}$ is $\bm{\phi}$ subject to a lower bound on the value of $\phi_{\sdot}$, implicitly set by all other components of $\bm{\phi}$---if initially $f_{\sdot}\left(\bmm{t}^{\circledast}\right) < 0$, then $\phi_{\sdot}$ must be effectively increased until $f_{\sdot}\left(\bmm{t}^{\circledast}\right) = 0$~\cite{molesky2020hierarchical}. 
That is, to evaluate $\mathcal{D}\left(\bm{\phi}\right)$ at any $\bm{\phi}\in\Phi_{\epsilon}$, set $\bmm{t}^{\circledast} = \bmm{A}_{\bm{\phi}}^{-1}\textbf{y}_{\bm{\phi}}$ (using the pseudo inverse in the event that $\bmm{A}_{\bm{\phi}}$ is semi-definite), and then check if $f_{\sdot}\left(\bmm{t}^{\circledast}\right) \geq 0$. 
If this condition is satisfied, or $\bmm{A}_{\bm{\phi}}$ is semi-definite, then Eq.~\eqref{dualEval} can be used directly. 
If not, then the value of $\mathcal{D}$ can be computed by increasing the effective value of $\phi_{\sdot}$ used in $\bmm{t}^{\circledast} = \bmm{A}_{\tilde{\bm{\phi}}}^{-1}\textbf{y}_{\tilde{\bm{\phi}}}$ (but not the actual value of $\phi_{\sdot}$) until $f_{\sdot}\left(\bmm{t}^{\circledast}\right) = 0$~\cite{molesky2020hierarchical}.
The resulting $\bmm{t}^{\circledast}$ and $\bmm{A}_{\tilde{\bm{\phi}}}$ can then be used in Eq.~\eqref{dualEval}.
Again, the implicit dependence of $\phi_{\sdot}$ is only required as a means of evaluation, and is not part of the definition of $\Phi_{\epsilon}$.
Moreover, in practice, is not actually necessary to treat $\phi_{\sdot}$ as being implicitly determined since the mechanics of the minimization will guarantee that $f_{\sdot}\left(\bmm{t}^{\circledast}\right) \geq 0$ at termination---the same dual minimum is determined regardless of whether or not one enforces that $\bmm{t}^{\circledast}\in\mathtt{C}$ by the corollary of lemma 3. 
If $\epsilon > 0$, Eq.~\eqref{dualEval} is smooth~\cite{zualinescu2021canonical} with
\begin{align}
  &\frac{d\mathcal{D}\left(\bm{\phi}\right)}{d\phi_{k}} = f_{k}\left(\bmm{t}^{\circledast}\right),
  \label{dualDers} \\
  &\frac{d^{2}\mathcal{D}\left(\bm{\phi}\right)}{d\phi_{k}d\phi_{j}} = 
  2\Re\left[\left(\bmm{s}_{k}-\bmm{A}_{k}\bmm{t}^{\circledast}\right)^{\dagger}
  \bmm{A}_{\bm{\phi}}^{-1}
  \left(\bmm{s}_{j}-\bmm{A}_{j}\bmm{t}^{\circledast}\right)\right].
  \nonumber
\end{align}
The following lemma gives a sufficient condition for $\min_{\bm{\phi}\in\Phi_{\epsilon}}\mathcal{D}\left(\bm{\phi}\right)$ to be defined, implying that $\inf$ can be replaced by $\min$ for optimization of the dual.
\\ \\
\emph{Lemma 2}---If $\exists\delta > 0$ such that for every vector of unit magnitude in the subset $\mathbb{S}^{\tilde{\mathtt{K}}}\cap\left\{\bm{\psi}/\lVert\bm{\psi}\rVert_{_\mathtt{E}}~\Big|~\bm{\psi}\in\Phi_{\epsilon}\right\}$
\begin{align}
 &\left[\max_{\bmm{t}\in\mathtt{C}}~
 2~\Re\left(\bmm{t}^{\dagger}\tilde{\bmm{s}}_{\hat{\bm{\psi}}}
 \right) - \bmm{t}^{\dagger}\tilde{\bmm{A}}_{\hat{\bm{\psi}}}\bmm{t}\right] > \delta,
 \nonumber
\end{align}
where $\tilde{\bmm{s}}_{\hat{\bm{\psi}}} = \sum_{k\in\tilde{\mathtt{K}}} \hat{\psi}_{k}\bmm{s}_{k}$, and $\tilde{\bmm{A}}_{\hat{\bm{\psi}}} = \sum_{k\in\tilde{\mathtt{K}}} \hat{\psi}_{k}\bmm{A}_{k}$, there is a $\bm{\phi}\in\Phi_{\epsilon}$ achieving the infimum of $\mathcal{D}\left(\bm{\phi}\right)$. 
\begin{proof}[Proof of lemma 2.]
  Because $\mathtt{C}$ is compact, for any optimization of the form of Eq.~\eqref{waveQCQP}, there is always finite number $P = \max_{\bmm{t}\in\mathtt{C}}f_{\circ}\left(\bmm{t}\right)$ bounding $\inf_{\bm{\phi}\in\Phi_{\epsilon}}\mathcal{D}\left(\bm{\phi}\right)$ from above, and a finite number $N = \min_{\bm{t}\in\mathtt{C}}f_{\circ}\left(\bmm{t}\right)$ bounding $f_{\circ}$ from below. 
  Take $\mathtt{B}_{\rho}$ to be a closed ball in $\mathbb{R}^{\tilde{\mathtt{K}}}$ centered on the origin of radius $\rho$. 
  Set $\rho = 1 + \left(P + N\right)/\delta$, and suppose that $m\hat{\bm{\psi}}\in\Phi_{\epsilon}$ and $m\hat{\bm{\psi}}\not\in \mathtt{B}_{\rho}\cap\Phi_{\epsilon}$.
  Then, 
  \begin{align}
   \mathcal{D}\left(m\hat{\bm{\psi}}\right) &=
   \max_{\bmm{t}\in\mathtt{C}}~m\left[2~\Re\left(\bmm{t}^{\dagger}\bmm{s}_{\hat{\bm{\psi}}}\right) - 
   \bmm{t}^{\dagger}\bmm{A}_{\hat{\bmm{\psi}}}
   \bmm{t}+\frac{1}{m}f_{\circ}\left(\bmm{t}\right)\right]
   \nonumber \\
   &\geq \rho\delta - N = P + \delta. 
   \nonumber
  \end{align}
  Therefore, as $\Phi_{\epsilon}$ is closed, $\left\{\bm{\phi}\in\Phi_{\epsilon}~|~\mathcal{D}\left(\bm{\phi}\right)\leq P\right\}$ is nonempty and bounded, and $\mathcal{D}$ is a continuous function of $\bm{\phi}$, the lemma is verified~\cite{bredon2013topology}. 
\end{proof}
Take $\mathtt{S}$ to be the collection of all vectors defining the linear part of a quadratic constraint contained in $\tilde{\mathtt{K}}$. 
The conditions of lemma 2 are then satisfied whenever $\mathtt{S}$ is linearly independent, and the additional restriction specified by the set $\left\{\bm{\psi}/\lVert\bm{\psi}\rVert_{_\mathtt{E}}~\Big|~\bm{\psi}\in\Phi_{\epsilon}\right\}$ limits admissible linear combinations $\bmm{s}_{\bm{\psi}}$ in a way such that there is always as $\bmm{c}\in\mathtt{C}$ with $\Re\left(\bmm{c}^{\dagger}\bmm{s}\right) > 0$. 
The assumption that a dual minimizer exists will often be used hereafter. 
\\ \\
\emph{Lemma 3}---Assuming that a dual minimizer exists, whenever $\epsilon > 0$ the $\bmm{t}$ determined by $\min_{\bm{\phi}\in\Phi_{\epsilon}}\max_{\bmm{t}\in\mathtt{C}}\mathcal{L}\left(\bm{\phi},\bmm{t}\right)$ is identical to a $\bmm{t}$ determined by $\max_{\bmm{t}\in\mathtt{C}}\inf_{\bm{\phi}\in\Phi_{\epsilon}}\mathcal{L}\left(\bm{\phi},\bmm{t}\right)$. 
For this $\bmm{t}$, the dual minimizer also achieves the infimum of $\max_{\bmm{t}\in\mathtt{C}}\inf_{\bm{\phi}\in\Phi_{\epsilon}}\mathcal{L}\left(\bm{\phi},\bmm{t}\right)$
\\ \\
\emph{Corollary}---If $\bm{\phi}^{\circledast}$ is a minimizer of $\inf_{\bm{\phi}\in\Phi_{\epsilon}}\mathcal{D}\left(\bm{\phi}\right)$, and $\bmm{t}^{\circledast}$ is a solution of $\max_{\bmm{t}\in\mathtt{C}}\mathcal{L}\left(\bm{\phi}^{\circledast},\bmm{t}\right)$, then $f_{\bm{\psi}}\left(\bmm{t}^{\circledast}\right) \geq 0$ for every positive semi-definite constraint $f_{\bm{\psi}}$ that can be formed by combining the constraints of $\tilde{\mathtt{K}}$.
\begin{proof}[Proof of lemma 3.]
  Let $\tilde{\bm{\phi}}^{\circledast}$ be the minimizer of $\mathcal{D}\left(\bm{\phi}\right)$ over $\Phi_{\epsilon}$. 
  Take $\tilde{\bmm{t}}^{\circledast}$ to be the corresponding, uniquely determined, $\bmm{t}\in\mathtt{C}$ as described in the discussion proceeding Eq.~\eqref{dualEval}. 
  By Sion's minimax theorem, there is then a $\bmm{t}^{\circledast}\in\mathtt{C}$ such that 
  $
    \inf_{\bm{\phi}\in\Phi_{\epsilon}}\mathcal{L}\left(\bm{\phi},\bmm{t}^{\circledast}\right) = 
    \mathcal{L}(\tilde{\bm{\phi}}^{\circledast},\tilde{\bmm{t}}^{\circledast}). 
  $
  Because $\tilde{\bmm{t}}^{\circledast}$ is selected by $\max_{\bmm{t}\in\mathtt{C}}\mathcal{L}\left(\tilde{\bm{\phi}}^{\circledast},\bmm{t}\right)$, 
  $$
   \mathcal{L}\left(\tilde{\bm{\phi}}^{\circledast},\tilde{\bmm{t}}^{\circledast}\right) = \inf_{\bm{\phi}\in\Phi_{\epsilon}}\mathcal{L}\left(\bm{\phi},\bmm{t}^{\circledast}\right)\leq 
   \mathcal{L}(\tilde{\bm{\phi}}^{\circledast},\bmm{t}^{\circledast}) \leq 
   \mathcal{L}\left(\tilde{\bm{\phi}}^{\circledast},\tilde{\bmm{t}}^{\circledast}\right).
  $$
  Therefore, $\mathcal{L}(\tilde{\bm{\phi}}^{\circledast},\bmm{t}^{\circledast}) = \mathcal{L}(\tilde{\bm{\phi}}^{\circledast},\tilde{\bmm{t}}^{\circledast})\Rightarrow\bmm{t}^{\circledast} = \tilde{\bmm{t}}^{\circledast}$, and $\mathcal{L}\left(\tilde{\bm{\phi}}^{\circledast},\tilde{\bmm{t}}^{\circledast}\right) = 
  \inf_{\bm{\phi}\in\Phi_{\epsilon}}\mathcal{L}\left(\bm{\phi},\tilde{\bmm{t}}^{\circledast}\right)$.
\end{proof} 
\subsection*{Condition for Strong Duality in Compact QCQPs}
The main insight offered by Sion's minimax theorem, in the present setting, is that whenever a positive definite quadratic equality constraint remains persistently satisfied under the relaxation of a compact QCQP to $\max_{\bmm{t}\in\mathtt{C}}\inf_{\bmm{\phi}\in\Phi_{\epsilon}}\mathcal{L}\left(\bm{\phi},\bmm{t}\right)$, then strong duality holds. 
In detail, assuming $\epsilon > 0$ and that a dual minimizer exists, suppose that the domain of possible maxima for $\inf_{\bm{\phi}\in\Phi_{\epsilon}}\mathcal{L}\left(\bm{\phi},\bmm{t}\right)$ can be restricted to a subset of $\mathtt{C}$ such that for every $\bmm{t}$ in the subset there exists an $\epsilon$-positive definite quadratic constraint, $f_{d} = \psi_{\sdot}^{d}f_{\sdot}+\sum_{_{\mathtt{K}}}\psi_{k}^{d}f_{k}$ for some $\bm{\psi}^{d}$, such that $f_{d}\left(\bmm{t}\right) = 0$. 
If at any such $\bmm{t}$ there is a constraint $f_{k}\in\tilde{\mathtt{K}}$ such that $f_{k}\left(\bmm{t}\right) \neq 0$, then there exists a scalar $d_{k}$ such that $\pm f_{k}\left(\bmm{t}\right) + d_{k}f_{d}\left(\bmm{t}\right) < 0$ and $\mp \bmm{A}_{k} + d_{k}\bmm{A}_{\bm{\psi}_{d}}\succeq \epsilon$. 
Referring to the combination $\pm f_{k} + d_{k}f_{d}$ as $f_{d'}$, by increasing $\bm{\phi}$ along the multiplier direction specified by $\bm{\phi}_{d'}$ the value of $\mathcal{L}\left(\bm{\phi},\bmm{t}\right)$ can then be made arbitrarily negative within $\Phi_{\epsilon}$, meaning that $\inf_{\bm{\phi}\in\Phi}\mathcal{L}\left(\bm{\phi},\bmm{t}\right) = -\infty$. 
Because feasible points are assumed to exist, and at feasible point $f_{\circ}$ must take on a finite value, $\max_{\bmm{t}\in\mathtt{C}}\mathcal{F}_{\epsilon}\left(\bmm{t}\right)$, with 
\begin{equation}
  \mathcal{F}_{\epsilon}\left(\bmm{t}\right) = 
  \inf_{\bm{\phi}\in\Phi_{\epsilon}}\mathcal{L}\left(\bm{\phi},\bmm{t}\right),
  \label{revDualDef}
\end{equation} 
will occur at a feasible point. 
Lemma 3 ensures that dual the solution is then also feasible, proving strong duality. 
\subsection*{``Dual'' Properties under Minimax Exchange}
Independent of the above condition for strong duality, the function $\mathcal{F}_{\epsilon}\left(\bmm{t}\right)$ is of clear theoretical value for characterizing the behaviour of duality transformations for compact QCQPs.
Decomposing $\mathcal{F}_{\epsilon}$ as 
\begin{align}
  \mathcal{F}_{\epsilon}\left(\bmm{t}\right) &= f_{\circ}\left(\bmm{t}\right) + q_{\epsilon}\left(\bmm{t}\right),
  \\
  q_{\epsilon}\left(\bmm{t}\right) &=\inf_{\bm{\phi}\in\Phi_{\epsilon}}\mathcal{L}\left(\bm{\phi},\bmm{t}\right) - f_{\circ}\left(\bmm{t}\right),
  \nonumber
\end{align}
define 
\begin{align}
  &\mathtt{Q}= \left\{\bmm{t}\in\mathtt{C}~\Big|~q_{\epsilon}\left(\bmm{t}\right)\geq 0\right\},
  \label{sionSets}\\
  &\mathtt{Q}_{_{0}} =\left\{\bmm{t}\in\mathtt{C}~\Big|~q_{\epsilon}\left(\bmm{t}\right) = 0\right\} \subseteq \mathtt{Q}.
\end{align}
Making use of the preceding reasoning, strong duality in a compact QCQP necessarily follows if there is an $\epsilon > 0$ resulting in $\bmm{t}^{\circledast}\in\mathtt{Q}_{_{0}}$ under the assumption that $\bmm{A}_{\circ}\preceq 0$. 
Moreover, in general, the following properties hold. 
\\ \\
\emph{Lemma 4}---$q_{\epsilon}\left(\bmm{t}\right) - \bmm{t}^{\dagger}\bmm{A}_{\circ}\bmm{t}$ is a concave function of $\bmm{t}$. 
If $\bmm{A}_{\circ}\preceq 0$, then $q_{\epsilon}\left(\bmm{t}\right)$ is concave, $\mathtt{Q}$ is closed and convex, $\bmm{t}\in\mathtt{Q}_{0}$ implies that $\bmm{t}$ satisfies all imposed constraints, and $\bmm{t}\not\in\mathtt{Q}$ implies $q_{\epsilon}\left(\bmm{t}\right) = -\infty$. 
\\ \\
\begin{proof}[Proof of Lemma 4.]
 Take $\bmm{t}_{1},\bmm{t}_{2}\in\mathtt{C}$ and $\alpha\in\left(0,1\right)$. 
 For any $\bmm{t} = \left(1-\alpha\right)\bmm{t}_{1} +\alpha\bmm{t}_{2}$ and $\sigma > 0$, let $f_{\bm{\phi}}\left(\bmm{t}\right)$ be quadratic function such that $f_{\circ}\left(\bmm{t}\right) + f_{\bm{\phi}}\left(\bmm{t}\right)$ is $\epsilon$-positive definite and $\lVert f_{\bm{\phi}}\left(\bmm{t}\right) - q_{\epsilon}\left(\bmm{t}\right)\rVert_{_\mathtt{E}} \leq \sigma$. 
 Then, $\forall \sigma$
 \begin{align}
   &f_{\circ}\left(\bmm{t}\right) + q_{\epsilon}\left(\bmm{t}\right) +\sigma \geq f_{\circ}\left(\bmm{t}\right) + f_{\bm{\phi}}\left(\bmm{t}\right)  \geq 
   \nonumber \\
   &\left(1-\alpha\right)\left(f_{\circ}\left(\bmm{t}_{1}\right) + f_{\bm{\phi}}\left(\bmm{t}_{1}\right)\right) + 
   \alpha\left(f_{\circ}\left(\bmm{t}_{2}\right) + f_{\bm{\phi}}\left(\bmm{t}_{2}\right)\right)\Rightarrow
   \nonumber \\
   &f_{\circ}\left(\bmm{t}\right) + q_{\epsilon}\left(\bmm{t}\right) \geq 
   \nonumber \\
   &\left(1-\alpha\right)\left(f_{\circ}\left(\bmm{t}_{1}\right) + q_{\epsilon}\left(\bmm{t}_{1}\right)\right) + 
   \alpha\left(f_{\circ}\left(\bmm{t}_{2}\right) + q_{\epsilon}\left(\bmm{t}_{2}\right)\right) - \sigma
   \nonumber 
 \end{align} 
 Therefore, since $2\Re\left(\bmm{t}^{\dagger}\bmm{s}_{\circ}\right)$ is linear, $q_{\epsilon}\left(\bmm{t}\right) - \bmm{t}^{\dagger}\bmm{A}_{\circ}\bmm{t}$ is concave function of $\bmm{t}$. 
 If $\bmm{A}_{\circ}\succeq 0$, then $f_{\circ}$ is convex, and the continued validity of the inequalities stated above proves that $q_{\epsilon}$ is concave. 
 \\ \\
 Keeping the assumption that $\bmm{A}_{\circ}\preceq 0$, let 
 $$
    p_{\rho}^{\epsilon}\left(\bmm{t}\right) = \inf_{\bm{\phi}\in\Psi_{\rho}^{\epsilon}}f_{\bm{\phi}}\left(\bmm{t}\right)
 $$
 with $\Psi^{\epsilon}_{\rho} = \left\{\bm{\phi}\in \mathtt{B}_{\rho}~|~\bmm{A}_{\bm{\phi}}\succeq\epsilon\right\}$ and $\mathtt{B}_{\rho}$ a closed ball of radius $\rho > 0$ centered on the origin of $\mathbb{R}^{\tilde{\mathtt{K}}}$. 
 For any finite value of $\rho$, $p_{\rho}^{\epsilon}\left(\bmm{t}\right)$ is then properly defined, continuous, function for all $\bmm{t}\in\Omega$ and concave. 
 Hence, $\left(p_\rho^{\epsilon}\right)^{-1}\left[0,\infty\right)$ is convex and closed~\cite{rockafellar1970convex}, and, as 
 $$
  \mathtt{Q} = \bigcap_{\rho\in\mathbb{N}}\left(p_\rho^{\epsilon}\right)^{-1}\left[0,\infty\right),
 $$
 $\mathtt{Q}$ is also closed and convex. 
 If $\bmm{t}\in\mathtt{Q}_{0}$, then, for any $\delta > 0$, there must be an $\epsilon$-positive semi-definite constraint $0\leq f_{\bm{\phi}}\left(\bmm{t}\right) \leq \delta$. 
 Take any $f_{k}\in\tilde{\mathtt{K}}$ and set $\lambda_{k} = \lVert \bmm{A}_{k} \rVert_{_\mathtt{O}}$. 
 $\pm\left(\epsilon/\lambda_{k}\right) f_{k} + 2f_{\bm{\phi}}$ is then $\epsilon$-positive definite, implying that $\pm\left(\epsilon/\lambda_{k}\right) f_{k}\left(\bm{t}\right) + 2f_{d}\left(\bm{t}\right) \geq 0\Rightarrow$  
 \begin{equation}
    \left(\forall\delta >0\right)~\delta \frac{2\lambda_{k}}{\epsilon}\geq \lVert f_{k}\left(\bmm{t}\right)\rVert_{_\mathtt{E}}.
    \nonumber
 \end{equation}
 If $\bmm{t}\not\in\mathtt{Q}$, then there is an $\epsilon$-positive definite constraint that is negative. 
 Any positive multiple of this constraint is then also $\epsilon$-positive definite, and so $q_{\epsilon}\left(\bmm{t}\right) = -\infty$. 
\end{proof}
In regards to lemma 4 and its implications for strong duality, it is important to realize that generally $\mathtt{Q}_{0} \neq \mathtt{Q}^{\partial}$. 
Given that $\bmm{t}\in\mathtt{Q}_{0}$ implies that $\bmm{t}$ satisfies all imposed constraints, $\bmm{t}\in\mathtt{Q}_{0}\Rightarrow\bmm{t}\in\mathtt{C}^{\partial}$. 
However, $\mathtt{Q}$ is convex, and so, whenever $\mathtt{Q}\neq\mathtt{C}$ it must be that $\mathtt{Q}_{0}\neq\mathtt{Q}^{\partial}$. 
That is, if $\mathtt{Q}\neq\mathtt{C}$, then some portion of $\mathtt{Q}^{\partial}$ must lie within the interior of $\mathtt{C}$, meaning that $f_{\sdot}\left(\bmm{t}\right) \neq 0$ for some $\bmm{t}\in\mathtt{Q}^{\partial}$---the constraint specified by $f_{\sdot}$ is not satisfied on all of $\mathtt{Q}^{\partial}$, and so $\mathtt{Q}^{\partial}\neq\mathtt{Q}_{0}$. 
It should also be noted that by selecting a suitable value of $c_{\sdot} > 0$ the redefinition of the objective as $f_{\circ} - c_{\sdot}f_{\sdot}$ can always be employed in order to meet the conditions stated for the second part of the lemma. 
Because $f_{\sdot}$ is an admissible constraint in the infimum over $\Phi_{\epsilon}$, the $\max_{\bmm{t}\in\mathtt{C}}\mathcal{F}_{\epsilon}\left(\bmm{t}\right)$ formulation shows that such a change can only ever reduce the value of the dual, becoming trivial if strong duality occurs. 
Because $f_{\sdot}$ is a constraint of the primal problem, the altered problem statement remains a bound on the unaltered QCQP. 
This same reasoning can be applied to any $\epsilon$-positive definite constraint in $\tilde{\mathtt{K}}$.

Given these four lemmas, we believe that it may be reasonably conjectured that in many QCQPs motivated by physical device design strong duality should hold.  
In analogy with the scattering potential optimization problems described earlier---where the magnitude of the field generated by any passive potential $\bmm{V}_{\rho}$ in response to an incident field is bounded by the conservation of resistive (``real'') power---for almost any practical open system it is possible to impose a positive definite quadratic equality constraint so that the analysis of the past two sections is applicable.
This then suggests, more often than not, that at least one compact constraint will remain satisfied under the relaxation of the QCQP to $\max_{\bmm{t}\in\mathtt{C}}\mathcal{F}_{\epsilon}\left(\bmm{t}\right)$.
For instance, if $f_{\sdot}\left(\bmm{t}^{\circledast}\right) > 0$, then $\bmm{t}^{\circledast}$ exists within the interior $\mathtt{C}$ as a field with ``excessive resistive power''. 
Therefore, $\bmm{t}^{\circledast}$ would also be a feasible point for the QCQP $\max_{\bmm{t}\in\Omega} f_{\mathsf{o}}\left(\bmm{t}\right) \ni \tilde{f}_{\sdot}\left(\bmm{t}\right) = 0$, where $\tilde{f}_{\sdot}$ is $f_{\sdot}$ subject to an artificial increase in the rate of irrecoverable energy loss, e.g. as quantified by $i\bmm{V}_{d}^{-1\mathsf{s}}$ in the scattering problem. 
Empirically, at least in the scattering problems we are most familiar with~\cite{gustafsson2019maximum,molesky2020t,kuang2020maximal}, constraints related to the conservation of resistive power are strong restrictions on most objectives related to physical device performance---loosely, realizable response characteristics for practically interesting design problems are often limited by material absorption and coupling to the environment---and so, the possibility that $\bmm{t}^{\circledast}$ lies ``far'' from $\mathtt{C}^{\partial}$ is innately odd. 
Moreover, in the nontrivial case that more than a single constraint is used, a $\bmm{t}^{\circledast}$ ``deep'' within the interior of $\mathtt{C}$ is likely nearer to the boundary of some other $\epsilon$-positive definite constraint than it is to $\mathtt{C}^{\partial}$. 
Hence, by redefining the objective as suggested above, it is also well possible that such a point can be ruled out completely based on $\tilde{\mathtt{K}}$.
Conversely, if $\bmm{t}^{\circledast}$ is ``near'' $\mathtt{C}^{\partial}$, then there is an $\epsilon$-positive definite constraint that is almost satisfied.
This, in turn, imposes restrictions on the degree to which any constraint in $\tilde{\mathtt{K}}$ can be violated as shown in the proof of lemma 4: the ``advantage'' that may be gained from violating some of the imposed constraints becomes increasingly small as $\bmm{t}$ tends to $\mathtt{C}^{\partial}$. 
With no obvious benefit offered by alternative positions, the supposition that the $\bmm{t}$ determined by $\max_{\bmm{t}\in\mathtt{C}}\mathcal{F}_{\epsilon}\left(\bmm{t}\right)$ should in most cases occur in $\mathtt{Q}_{0}\cap \mathtt{C}^{\partial}$ under trivial changes to the objective is, at least to the authours, quite plausible. 
\subsection*{Attaining Strong Duality by Objective Modifications}
The perspective of the previous section recommends several means by which any compact, feasible, QCQP can be retroactively modified until it becomes strongly dual---if $\bmm{s}_{\mathsf{o}}$ is modified as detailed below then the optimal $\bmm{t}^{\circledast}$ of $\max_{\bmm{t}\in\mathtt{C}}\mathcal{F}_{\epsilon}\left(\bmm{t}\right)$ will occur in $\mathtt{Q}_{0}$ and strong duality will hold. 
The reasoning supporting this idea is based on two features of Eq.~\eqref{sion}. 
First, the value of $\bmm{s_{\circ}}$ does not enter into either $\Phi_{\epsilon}$ or $\bmm{q}_{\epsilon}$.  
$\bmm{s}_{\circ}$, correspondingly, can be viewed as an independent parameter for altering the value of either $\mathcal{F}_{\epsilon}$ (at a given $\bmm{t}$) or $\mathcal{D}$ (at a given $\bm{\phi}$). 
Second, it is clear from the $\max_{\bmm{t}\in\mathtt{C}}\mathcal{F}_{\epsilon}\left(\bmm{t}\right)$ perspective that if $\bmm{s}_{\circ}$ becomes sufficiently large and sufficiently aligned with any feasible point, then that feasible point will become the solution of $\min_{\bm{\phi}\in\Phi_{\epsilon}}\mathcal{D}\left(\bm{\phi}\right)$---by changing $\bmm{s}_{\circ}$, $\mathcal{F}_{\epsilon}\left(\bmm{t}\right)$ can be made to increase when moving towards some feasible point. 

Recalling that it is not necessary to treat $\phi_{b}$ as being implicitly defined for algorithmic considerations, since minimization will directly guarantee that at termination $f_{\sdot}\left(\bmm{t}^{\circledast}\right) \geq 0$, there are two possible outcomes for $\min_{\bm{\phi}\in\Phi_{\epsilon}}\mathcal{D}\left(\bm{\phi}\right)$: the convex solver terminates with $\bm{\phi}^{\circledast}$ such that (a) $\bmm{A}_{\bm{\phi}^{\circledast}}\succ\epsilon$, or (b) $\bmm{A}_{\bm{\phi}^{\circledast}} -\epsilon\bmm{I}$ is semi-definite. 

In the first case it is simple to establish that strong duality holds by Eq.~\eqref{dualDers}~\cite{zualinescu2021canonical,molesky2020hierarchical}. 
In the case second, knowing that feasible points exist, it possible to propose a variety of auxiliary optimization problems that either (a) alters the current $\bmm{t}^{\circledast}$ in such a way that dual minimization can be restarted with all constraint violations reduced, or (b) strong duality is realized outright. 
For example, take $v_{p}$ to denote the value of $f_{p}\in\mathtt{K}$, $v_{\sdot}$ to be the value of 
$f_{\sdot}$, and let $a_{\circ}$, $\delta_{p}$ and $\delta_{\sdot}$ be real numbers greater than zero; as an illustration of the first approach one may determine a feasible point of the QCQP    
\begin{align}
  &\max_{\bmm{t}_{\delta}\in\Omega}~2~\Re\left(\bmm{t}_{\delta}^{\dagger}\bmm{s}_{\circ}\right)-a_{\circ}~\bmm{t}^{\dagger}_{\delta}\bmm{t}_{\delta}\ni
  \\
  &\left(\forall f_{p}\in\mathtt{K}\right) -\left|v_{p}\right|+\delta_{p}\leq f_{p}\left(\bmm{t} +\bmm{t}_{\delta}\right) \leq \left|v_{p}\right|-\delta_{p}~\land
  \nonumber \\
  &f_{\sdot}\left(\bmm{t} +\bmm{t}_{\delta}\right) \leq v_{\sdot} - \delta_{\sdot}.
  \nonumber
\end{align}
If $\delta_{\sdot}$ is taken to be slightly larger than $v_{\sdot}$, then, by Eq.~\eqref{dualDers}, for a sufficiently small increase in $\phi_{\sdot}$, in the modified QCQP determined by $s_{\circ}\rightarrow s_{\circ} + \bmm{A}_{\bm{\phi}^{\circledast}_{l}}^{-1}\bmm{t}_{\delta}$, where $\bm{\phi}^{\circledast}_{l}$ denotes the value of $\bm{\phi}$ at the latest dual minimum, $\bmm{A}_{\bm{\phi}}\succ\epsilon$, and every imposed constraints will be less violated than at the last termination. 
(Here, the value of $a_{\circ}$ controls the relative importance of the modification being small compared to the relative importance that the modification is well suited to the original objective.) 
Therefore, dual minimization can be restarted.
As basic implementation of the more direct approach, one may determine a ``single shot'' $\bmm{t}^{\delta}$ such all constraints are satisfied by local gradient descent. 

Hybrids of these two methods, in conjunction with the subtraction of positive definite constraints from the objective as sketched above, are also possible. 
Specifically, one could begin by calculating a $\bmm{t}^{\delta}$, and an associated source modification $\bmm{s}^{\delta}$, that leads to the immediate satisfaction of all constraints. 
Then, rather than just using the modification directly, one could determine the smallest value of $\alpha_{\circ}\in\left(0,1\right]$ such that attempting to minimizing the dual over $\bm{\phi}\in\Phi_{\epsilon}$ leads to a change in $\bm{\phi}^{\circledast}$. 
\subsection*{QCQP Updates and Modified Strong Duality}
As a final remark, there is also a potentially interesting interaction between the original compact QCQP, and any modified version of the QCQP that is strongly dual. 
Namely, once the modified problem $\mathcal{Q}^{m}\left(\Omega,f_{\circ}^{m},\tilde{\mathtt{K}}\right)$ (with an $\!^{m}$ superscript denoting modification) has been solved, the value of $f_{\circ}^{m}$ within the feasible set determined by $\tilde{\mathtt{K}}$ is bounded by the value of the modified program: 
suppose $\left\{\bmm{t}\in \Omega~\big|~\left(\forall f_{k}\in\tilde{\mathtt{K}}\right)~f_{k}\left(\bmm{t}\right) =0\right\}$, then
\begin{equation}
  f^{m}_{\circ}\left(\bmm{t}\right)\leq 
  \mathcal{Q}^{m}\left(\Omega,f_{\circ}^{m},\tilde{\mathtt{K}}\right).
  \label{modQCQPInfo}
\end{equation}
Because $f_{\circ}^{m}$ remains a quadratic function, Eq.~\eqref{modQCQPInfo} can be introduced as an additional, generally nontrivial, inequality constraint on the original compact QCQP. 
The inclusion of this extra information in no way alters the validity for the points discussed in previous section, leading to improved bounds, and, if so desired, new modifications for achieving strong duality. 
An iterative strategy along these lines suggests itself as a promising approach for discovering nearly optimal feasible points.   
\subsection*{Outlook}
The interpretation of compact QCQPs offered by Sion's theorem paints an expressly positive picture for future applications. 
While there may well be compact QCQPs in which strong duality does not hold, the relation of such dual solutions in physical contexts to polarization fields that do not sufficiently consume resistive power indicates both that these instances are probably somewhat special, and that, in any event, the size of the duality gap is probably not that large. 
Moreover, the extent of allowed constraint violation is directly tied to the positivity of all possible positive definite quadratic constraints that can be formed from the imposed constraint set. 
When any one positive definite quadratic constraint of this set is nearly satisfied, then the possible violation in any imposed constraints becomes increasingly bounded.  
Beyond of these considerations, as described in the last two sections, there are a number of reliable ways to modify the optimization object of any compact, feasible, QCQP in order to achieve strong duality. 
To the authours, this understanding strongly motivates the conception of a new class of inverse design algorithms for wave scattering applications that co-calculates limits and designs---guaranteeing that discovered solutions are either optimal or optimal to within some percentage~\cite{angeris2019computational,chao2021physical}. 
\subsection*{Acknowledgments}
This work was supported by the National Science Foundation under the Emerging Frontiers in Research and Innovation (EFRI) program, EFMA-1640986, the Cornell Center for Materials Research (MRSEC) through award DMR-1719875, the Defense Advanced Research Projects Agency (DARPA) under agreements HR00112090011, HR00111820046 and HR0011047197, and the Canada First Research Excellence Fund via the Institut de Valorisation des Données (IVADO) collaboration. 
The opinions and findings expressed within are those of the authors, and should not be interpreted as being representative of any institution. 
The authors thank Guillermo Angeris for several useful discussion concerning QCQPs, Sion's minimax theorem, and computational complexity.
\bibliography{esaLibF}

\begin{thebibliography}{44}
\providecommand{\natexlab}[1]{#1}
\providecommand{\url}[1]{\texttt{#1}}
\expandafter\ifx\csname urlstyle\endcsname\relax
  \providecommand{\doi}[1]{doi: #1}\else
  \providecommand{\doi}{doi: \begingroup \urlstyle{rm}\Url}\fi

\bibitem[Landau and Lifshitz(1976)]{landau1976mechanics}
Lev~Davidovich Landau and Evgenii~Mikhailovich Lifshitz.
\newblock \emph{Mechanics: Volume 1}, volume 1,~\S 1.
\newblock Butterworth-Heinemann, 1976.

\bibitem[Lu and Vu{\v{c}}kovi{\'c}(2012)]{lu2012objective}
Jesse Lu and Jelena Vu{\v{c}}kovi{\'c}.
\newblock Objective-first design of high-efficiency, small-footprint couplers
  between arbitrary nanophotonic waveguide modes.
\newblock \emph{Optics express}, 20\penalty0 (7):\penalty0 7221--7236, 2012.

\bibitem[Miller(2013)]{miller2013complicated}
David~AB Miller.
\newblock How complicated must an optical component be?
\newblock \emph{JOSA A}, 30\penalty0 (2):\penalty0 238--251, 2013.

\bibitem[Roes et~al.(2012)Roes, Duarte, Hendrix, and
  Lomonova]{roes2012acoustic}
Maurice~GL Roes, Jorge~L Duarte, Marcel~AM Hendrix, and Elena~A Lomonova.
\newblock Acoustic energy transfer: A review.
\newblock \emph{IEEE Transactions on Industrial Electronics}, 60\penalty0
  (1):\penalty0 242--248, 2012.

\bibitem[Simon(2015)]{simon2015quantum}
Barry Simon.
\newblock \emph{Quantum mechanics for Hamiltonians defined as quadratic forms},
  volume~68.
\newblock Princeton University Press, 2015.

\bibitem[Ganahl et~al.(2017)Ganahl, Rinc{\'o}n, and
  Vidal]{ganahl2017continuous}
Martin Ganahl, Juli{\'a}n Rinc{\'o}n, and Guifre Vidal.
\newblock Continuous matrix product states for quantum fields: An energy
  minimization algorithm.
\newblock \emph{Physical Review Letters}, 118\penalty0 (22):\penalty0 220402,
  2017.

\bibitem[Angeris et~al.(2021{\natexlab{a}})Angeris, Vu{\v{c}}kovi{\'c}, and
  Boyd]{angeris2021convex}
Guillermo Angeris, Jelena Vu{\v{c}}kovi{\'c}, and Stephen Boyd.
\newblock Convex restrictions in physical design.
\newblock \emph{Scientific Reports}, 11\penalty0 (1):\penalty0 1--10,
  2021{\natexlab{a}}.
\newblock \doi{10.1038/s41598-021-92451-1}.

\bibitem[Chao et~al.(2021)Chao, Strekha, Defo, Molesky, and
  Rodriguez]{chao2021physical}
Pengning Chao, Benjamin Strekha, Rodrick~Kuate Defo, Sean Molesky, and
  Alejandro~W. Rodriguez.
\newblock Physical limits on electromagnetic response.
\newblock \emph{arXiv:2109.05667}, 2021.

\bibitem[Angeris et~al.(2019)Angeris, Vu{\v{c}}kovi{\'c}, and
  Boyd]{angeris2019computational}
Guillermo Angeris, Jelena Vu{\v{c}}kovi{\'c}, and Stephen~P Boyd.
\newblock Computational bounds for photonic design.
\newblock \emph{ACS Photonics}, 6\penalty0 (5):\penalty0 1232, 2019.
\newblock \doi{10.1021/acsphotonics.9b00154}.

\bibitem[Gustafsson et~al.(2020)Gustafsson, Schab, Jelinek, and
  Capek]{gustafsson2020upper}
Mats Gustafsson, Kurt Schab, Lukas Jelinek, and Miloslav Capek.
\newblock Upper bounds on absorption and scattering.
\newblock \emph{New Journal of Physics}, 22\penalty0 (7):\penalty0 073013,
  September 2020.
\newblock ISSN 1367-2630.
\newblock \doi{10.1088/1367-2630/ab83d3}.

\bibitem[Molesky et~al.(2020{\natexlab{a}})Molesky, Chao, Jin, and
  Rodriguez]{molesky2020t}
Sean Molesky, Pengning Chao, Weiliang Jin, and Alejandro~W Rodriguez.
\newblock Global $\mathbb{T}$ operator bounds on electromagnetic scattering:
  upper bounds on far-field cross sections.
\newblock \emph{Physical Review Research}, 2\penalty0 (3):\penalty0 033172,
  2020{\natexlab{a}}.

\bibitem[Angeris et~al.(2021{\natexlab{b}})Angeris, Vu{\v{c}}kovi{\'c}, and
  Boyd]{angeris2021heuristic}
Guillermo Angeris, Jelena Vu{\v{c}}kovi{\'c}, and Stephen Boyd.
\newblock Heuristic methods and performance bounds for photonic design.
\newblock \emph{Optics Express}, 29\penalty0 (2):\penalty0 2827--2854,
  2021{\natexlab{b}}.

\bibitem[Liska et~al.(2021)Liska, Jelinek, and Capek]{liska2021fundamental}
Jakub Liska, Lukas Jelinek, and Miloslav Capek.
\newblock Fundamental bounds to time-harmonic quadratic metrics in
  electromagnetism: Overview and implementation.
\newblock \emph{arXiv 2110.05312}, 2021.

\bibitem[Aaronson and Is(2003)]{aaronson2003p}
Scott Aaronson and P~Is.
\newblock Is p versus np formally independent?
\newblock \emph{Bulletin of the EATCS}, 81\penalty0 (109-136):\penalty0 70,
  2003.

\bibitem[Park and Boyd(2017)]{park2017general}
Jaehyun Park and Stephen Boyd.
\newblock General heuristics for nonconvex quadratically constrained quadratic
  programming.
\newblock \emph{arXiv:1703.07870}, 2017.

\bibitem[Jensen and Sigmund(2011)]{jensen2011topology}
Jakob~S{\o}ndergaard Jensen and Ole Sigmund.
\newblock Topology optimization for nano-photonics.
\newblock \emph{Laser \& Photonics Reviews}, 5\penalty0 (2):\penalty0 308--321,
  2011.

\bibitem[Molesky et~al.(2018)Molesky, Lin, Piggott, Jin, Vu{\v{c}}kovi{\'c},
  and Rodriguez]{molesky2018inverse}
Sean Molesky, Zin Lin, Alexander~Y Piggott, Weiliang Jin, Jelena
  Vu{\v{c}}kovi{\'c}, and Alejandro~W Rodriguez.
\newblock Inverse design in nanophotonics.
\newblock \emph{Nature Photonics}, 12\penalty0 (11):\penalty0 659--670, 2018.

\bibitem[Christiansen and Sigmund(2021)]{christiansen2021inverse}
Rasmus~E Christiansen and Ole Sigmund.
\newblock Inverse design in photonics by topology optimization: tutorial.
\newblock \emph{JOSA B}, 38\penalty0 (2):\penalty0 496--509, 2021.

\bibitem[Digani et~al.(2022)Digani, Hon, and Davoyan]{digani2022framework}
Jagrit Digani, Philip~WC Hon, and Artur~R Davoyan.
\newblock Framework for expediting discovery of optimal solutions with blackbox
  algorithms in non-topology photonic inverse design.
\newblock \emph{ACS Photonics}, 9\penalty0 (2):\penalty0 432–442, 2022.
\newblock \doi{10.1021/acsphotonics.1c01819}.

\bibitem[Zhao et~al.(2022)Zhao, Boutami, and Fan]{zhao2022efficient}
Nathan~Z Zhao, Salim Boutami, and Shanhui Fan.
\newblock Efficient method for accelerating line searches in adjoint
  optimization of photonic devices by combining schur complement domain
  decomposition and born series expansions.
\newblock \emph{Optics Express}, 30\penalty0 (4):\penalty0 6413--6424, 2022.

\bibitem[Miller(2019)]{miller2019waves}
David~AB Miller.
\newblock Waves, modes, communications, and optics: a tutorial.
\newblock \emph{Advances in Optics and Photonics}, 11\penalty0 (3):\penalty0
  679--825, 2019.

\bibitem[Kuang et~al.(2020)Kuang, Zhang, and Miller]{kuang2020maximal}
Zeyu Kuang, Lang Zhang, and Owen~D. Miller.
\newblock Maximal single-frequency electromagnetic response.
\newblock \emph{Optica}, 7\penalty0 (12):\penalty0 1746--1757, Dec 2020.
\newblock \doi{10.1364/OPTICA.398715}.

\bibitem[Schab et~al.(2020)Schab, Rothschild, Nguyen, Capek, Jelinek, and
  Gustafsson]{schab2020trade}
Kurt Schab, Austin Rothschild, Kristi Nguyen, Miloslav Capek, Lukas Jelinek,
  and Mats Gustafsson.
\newblock Trade-offs in absorption and scattering by nanophotonic structures.
\newblock \emph{Optics Express}, 28\penalty0 (24):\penalty0 36584--36599, 2020.

\bibitem[Molesky et~al.(2020{\natexlab{b}})Molesky, Chao, and
  Rodriguez]{molesky2020hierarchical}
Sean Molesky, Pengning Chao, and Alejandro~W. Rodriguez.
\newblock Hierarchical mean-field $\mathbb{T}$ operator bounds on
  electromagnetic scattering: Upper bounds on near-field radiative purcell
  enhancement.
\newblock \emph{Physical Review Research}, 2:\penalty0 043398, Dec
  2020{\natexlab{b}}.
\newblock \doi{10.1103/PhysRevResearch.2.043398}.

\bibitem[Kuang and Miller(2020)]{kuang2020computational}
Zeyu Kuang and Owen~D. Miller.
\newblock Computational bounds to light--matter interactions via local
  conservation laws.
\newblock \emph{Physical Review Letters}, 125:\penalty0 263607, Dec 2020.
\newblock \doi{10.1103/PhysRevLett.125.263607}.

\bibitem[Jelinek et~al.()Jelinek, Gustafsson, Capek, and Schab]{jelinek2020sub}
L~Jelinek, M~Gustafsson, M~Capek, and K~Schab.
\newblock Sub-structure limits to optical phenomena.
\newblock In \emph{2020 Fourteenth International Congress on Artificial
  Materials for Novel Wave Phenomena (Metamaterials)}, pages 400--402. IEEE.

\bibitem[Capek et~al.(2021)Capek, Jelinek, and Masek]{capek2021fundamental}
Miloslav Capek, Lukas Jelinek, and Michal Masek.
\newblock Fundamental bounds for multi-port antennas.
\newblock In \emph{2021 15th European Conference on Antennas and Propagation
  (EuCAP)}, pages 1--4. IEEE, 2021.

\bibitem[Zhang et~al.(2021)Zhang, Kuang, Puri, and
  Miller]{zhang2021conservationlawbased}
Hanwen Zhang, Zeyu Kuang, Shruti Puri, and Owen~D. Miller.
\newblock Conservation-law-based global bounds to quantum optimal control.
\newblock \emph{Phys. Rev. Lett.}, 127:\penalty0 110506, Sep 2021.
\newblock \doi{10.1103/PhysRevLett.127.110506}.
\newblock URL \url{https://link.aps.org/doi/10.1103/PhysRevLett.127.110506}.

\bibitem[Taylor and Drori(2021)]{taylor2021optimal}
Adrien Taylor and Yoel Drori.
\newblock An optimal gradient method for smooth (possibly strongly) convex
  minimization.
\newblock \emph{arXiv:2101.09741}, 2021.

\bibitem[Tsang et~al.(2004)Tsang, Kong, and Ding]{tsang2004scattering}
Leung Tsang, Jin~Au Kong, and Kung-Hau Ding.
\newblock \emph{Scattering of electromagnetic waves: theories and
  applications}, volume~27.
\newblock John Wiley \& Sons, 2004.

\bibitem[Trivedi et~al.(2020)Trivedi, Angeris, Su, Boyd, Fan, and Vu\ifmmode
  \check{c}\else \v{c}\fi{}kovi\ifmmode~\acute{c}\else
  \'{c}\fi{}]{trivedi2020bounds}
Rahul Trivedi, Guillermo Angeris, Logan Su, Stephen Boyd, Shanhui Fan, and
  Jelena Vu\ifmmode \check{c}\else \v{c}\fi{}kovi\ifmmode~\acute{c}\else
  \'{c}\fi{}.
\newblock Bounds for scattering from absorptionless electromagnetic structures.
\newblock \emph{Physical Review Applied}, 14:\penalty0 014025, Jul 2020.
\newblock \doi{10.1103/PhysRevApplied.14.014025}.

\bibitem[Landau and Lifshitz(2013)]{landau2013statistical}
Lev~D Landau and Evgeny~M Lifshitz.
\newblock \emph{Statistical physics: volume 5}, volume~5.
\newblock Elsevier, 2013.

\bibitem[Felsen and Marcuvitz(1994)]{felsen1994radiation}
Leopold~B Felsen and Nathan Marcuvitz.
\newblock \emph{Radiation and scattering of waves}, volume~31.
\newblock John Wiley \& Sons, 1994.

\bibitem[Newton(2013)]{newton2013scattering}
Roger~G Newton.
\newblock \emph{Scattering theory of waves and particles}.
\newblock Springer Science \& Business Media, 2013.

\bibitem[Molesky et~al.(2022)Molesky, Chao, Mohajan, Reinhart, Chi, and
  Rodriguez]{molesky2021comm}
S.~Molesky, P.~Chao, J.~Mohajan, W.~Reinhart, H.~Chi, and A.~W. Rodriguez.
\newblock $\mathbb{T}$-operator limits on optical communication: Metaoptics,
  computation, and input-output transformations.
\newblock \emph{Physical Review Research}, 4:\penalty0 013020, Jan 2022.
\newblock \doi{10.1103/PhysRevResearch.4.013020}.

\bibitem[Bredon(2013)]{bredon2013topology}
Glen~E Bredon.
\newblock \emph{Topology and geometry}, volume 139.
\newblock Springer Science \& Business Media, 2013.

\bibitem[Boyd and Vandenberghe(2004)]{boyd2004convex}
Stephen Boyd and Lieven Vandenberghe.
\newblock \emph{Convex optimization}.
\newblock Cambridge University Press, 2004.

\bibitem[Sion(1958)]{sion1958general}
Maurice Sion.
\newblock On general minimax theorems.
\newblock \emph{Pacific Journal of Mathematics}, 8\penalty0 (1):\penalty0
  171--176, 1958.

\bibitem[Komiya(1988)]{komiya1988elementary}
Hidetoshi Komiya.
\newblock Elementary proof for sion's minimax theorem.
\newblock \emph{Kodai Mathematical Journal}, 11\penalty0 (1):\penalty0 5--7,
  1988.

\bibitem[Folland(1999)]{folland1999real}
Gerald~B Folland.
\newblock \emph{Real analysis: modern techniques and their applications},
  volume~40.
\newblock John Wiley \& Sons, 1999.

\bibitem[Rockafellar(1970)]{rockafellar1970convex}
R~Tyrrell Rockafellar.
\newblock \emph{Convex analysis}, volume~36.
\newblock Princeton University Press, 1970.

\bibitem[Beck and Eldar(2006)]{beck2006strong}
Amir Beck and Yonina~C Eldar.
\newblock Strong duality in nonconvex quadratic optimization with two quadratic
  constraints.
\newblock \emph{SIAM Journal on Optimization}, 17\penalty0 (3):\penalty0
  844--860, 2006.

\bibitem[Z{\u{a}}linescu(2021)]{zualinescu2021canonical}
Constantin Z{\u{a}}linescu.
\newblock On canonical duality theory and constrained optimization problems.
\newblock \emph{Journal of Global Optimization}, pages 1--18, 2021.

\bibitem[Gustafsson and Capek(2019)]{gustafsson2019maximum}
Mats Gustafsson and Miloslav Capek.
\newblock Maximum gain, effective area, and directivity.
\newblock \emph{IEEE Transactions on Antennas and Propagation}, 67\penalty0
  (8):\penalty0 5282, 2019.
\newblock \doi{10.1109/TAP.2019.2916760}.

\end{thebibliography}
\end{document}